\pdfoutput=1
\documentclass[conference]{IEEEtran}
\IEEEoverridecommandlockouts

\usepackage{eso-pic}
\usepackage{hyperref}
\usepackage{cite}
\usepackage{amsmath,amssymb,amsfonts}
\usepackage{graphicx}
\usepackage{textcomp}
\usepackage{xcolor}
\usepackage{amsthm}
\usepackage{url}
\usepackage{asymptote}
\hypersetup{hidelinks}

\newtheorem{thm}{Theorem}
\newtheorem{lemma}[thm]{Lemma}
\newtheorem*{mainthm}{Main Theorem}
\theoremstyle{definition}

\newcommand{\F}{\mathbb{F}_3}

\author{\IEEEauthorblockN{Justin Stevens}
\IEEEauthorblockA{\textit{Independent Researcher} \\
San Diego, CA, USA \\
justin@numbertheoryguy.com}}

\begin{document}

\title{The Maximum Number of Sets for 16 Cards is 26}

\maketitle

\begin{abstract}
In SET\textregistered{}, each card has four properties with three values each, and three cards form a \emph{set} when every property is all the same or all different. We prove that the maximum number of sets among $13$, $14$, $15$, and $16$ cards is $16$, $19$, $23$, and $26$, respectively. Viewing cards as points of $\mathbb{F}_3^4$ and sets as lines, we show that a configuration that leaves one of the $40$ line directions of $\mathbb{F}_3^4$ unused projects into $\mathbb{F}_3^3$ without losing a set, which reduces the problem to three properties, where exhaustive search is feasible. We force an unused direction by counting repeated directions in the planes through a card on many sets, and treat configurations containing a magic square separately.
\end{abstract}


\section{Introduction}

SET\textregistered{} is a card game with a deck of $81$ cards. In the standard version of the game, each card has four properties, each of which is capable of taking on three distinct values. A \textbf{set} comprises three cards in which each property is \textit{all the same or all different}. The game is played with $12$ cards to create a \textbf{board}. If no set is visible, the dealer adds three more cards, giving a board of $15$ cards.

We study the following question: \emph{what is the maximum number of sets that can appear among $n$ cards?} Write $M(n,k)$ for this maximum when cards have $k$ properties, so the standard game concerns $M(n,4)$. In earlier work with Wilson and Kelley~\cite{stevens2025maximum}, we gave elementary proofs of $M(n,4)$ for $3\le n\le 12$; in particular, $12$ cards contain at most $14$ sets, a fact first proven by de Bruijn~\cite{deBruijn2002Set}. A complete computer search with $4$ properties is impractical: for $n=12$, it means checking $\binom{81}{12}\binom{12}{3}\approx 1.5\cdot 10^{16}$ possible sets. In this paper, we extend the result to $16$ cards using elementary techniques. 

\begin{mainthm}
For every $3\le n\le 16$, the maximum number of sets among $n$ cards with four properties equals that with three properties: $M(n,4)=M(n,3)$, as follows:
\begin{center}\upshape
\setlength{\tabcolsep}{2.9pt}
\begin{tabular}{c|cccccccccccccc}
$n$ & 3 & 4 & 5 & 6 & 7 & 8 & 9 & 10 & 11 & 12 & 13 & 14 & 15 & 16\\ \hline
$M(n,4)$ & 1 & 1 & 2 & 3 & 5 & 8 & 12 & 12 & 13 & 14 & 16 & 19 & 23 & 26
\end{tabular}
\end{center}
\end{mainthm}

Section~\ref{sec:projection} reduces four properties to three when a direction is unused, and Section~\ref{sec:main} treats $13\le n\le 16$. The three-property maxima $M(n,3)$ were computed by exhaustive search in~\cite{stevens2025maximum}, with code on GitHub~\cite{setsearch}, and appear in OEIS A380167~\cite{oeis_A380167}. We establish that these values also give the four-property maxima for $3\le n\le 16$.

\section{Related Work}

We refer to~\cite{stevens2025maximum} for related work on SET\textregistered{}. Vinci~\cite{vinci2009maximum} found $16$, $19$, $23$, and $26$ sets for $13$, $14$, $15$, and $16$ cards by consecutive maximization, which gives only lower bounds. To our knowledge, this is the first proof that these values are the maximum in the four-property game.

\section{Preliminaries}

\subsection{Cards as points}

Let all boards with $n$ cards and $k$ properties (with $3$ values per property) be $\mathcal{B}(n,k)$. We can choose an arbitrary order for the properties, and for each property we can assign the numerical values $0$, $1$, and $2$ to its three values.  Given such an assignment, there is a bijection between the elements of $\mathcal{B}(n,k)$ and the sets $\{c_1, c_2, \cdots, c_n\}$ of $n$ distinct points $c_i\in\F^k$. Observe that in proofs throughout the rest of the paper, we refer to the cards on a board as points in $\F^k$.

Since the cards in a set must have each property all the same or all different, a set is a line in this vector space.

\begin{lemma}\label{thm:line} If we have two points $p$ and $q$, for a third point $r$ to form a set with them, it must satisfy: $p+q+r=0$ (in $\F^k$). \end{lemma}
\begin{proof}
We prove for each coordinate that the sum of the elements in $\F$ is $0$, so the sum is $0$ in $\F^k$. Observe that if all the values are the same, they must be $(0,0,0)$, $(1,1,1)$, or $(2,2,2)$. In all three cases, the sum of the elements is $0$ in $\F$. If all three values are different, then their numerical values are $(0,1,2)$ (or a permutation thereof), which sum to $0$ in $\F$. In both cases, the sum is $0$, so $p+q+r=0$.
\end{proof}

Conversely, three distinct points with $a+b+c=0$ form a set: in each coordinate, either all three values are different, or two values equal some $x$ and then the third equals $-2x=x$, so the values are all the same or all different. 

For a board $S$, let $N(S)$ be the number of sets in $S$, and let $M(n,k)$ be the maximum of $N(S)$ over all $S\in\mathcal{B}(n,k)$. We call $S$ a \emph{maximum} \emph{configuration} of $\mathcal{B}(n,k)$ if $N(S)=M(n,k)$. We use the values of $M(n,3)$ from the computer search of~\cite{stevens2025maximum}; for $3\le n\le 16$, they are the values listed in the Main Theorem. Note that $M(n,3)\le M(n,4)$, since we can set one property to be constant to get sets with $3$ properties.

Two distinct points $a$ and $b$ lie on exactly one line $\{a,\,a+v,\,a+2v\}$, where $v=b-a$, and we say that this line has \emph{direction} $v$. Its third point is $a+2v=-a-b$, so by Lemma~\ref{thm:line} and its converse, the sets of $S$ are exactly the lines contained in $S$. There are $(3^4-1)/2=40$ possible directions in $\F^4$, since the nonzero vectors $v$ and $2v$ represent the same direction. A direction $v$ is \emph{used} by $S$ if two points of $S$ differ by a multiple of $v$, and \emph{unused} otherwise.

\subsection{Magic squares}

\begin{figure}[t]
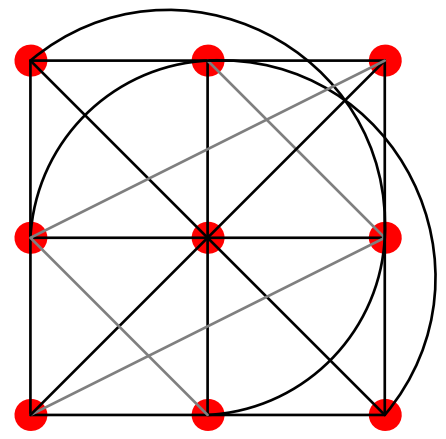

\centering
  \begin{asy}
    import graph;
    size(1.1inch);
    dot((0,0), 6+red);
    dot((1,0), 6+red);
    dot((0, 1), 6+red);
    dot((0,2), 6+red);
    dot((1, 1), 6+red);
    dot((1,2), 6+red);
    dot((2,1), 6+red);
    dot((2,0), 6+red);
    dot((2,2), 6+red);
    draw((0,0)--(0,1)--(0,2), black);
    draw((1,0)--(1,1)--(1,2), black);
    draw((0,1)--(1,1)--(2,1), black);
    draw((0,0)--(1,0)--(2,0), black);
    draw((0,2)--(1,1)--(2,0), black);
    draw((2,0)..(1,2)..(0,1), black);
    draw((0,0)--(2,1)--(1,2), grey);
    draw((1,0)..(2,1)..(0,2), black);
    draw((2,0)--(2,1)--(2,2), black);
    draw((0,2)--(1,2)--(2,2), black);
    draw((0,0)--(1,1)--(2,2), black);
    draw((1,0)--(0,1)--(2,2), grey);
  \end{asy}
  \caption{A magic square: $9$ cards/points that yield $12$ sets. Each line or curve forms a set, as do the two grey segments.}
  \label{fig:magicsquare}
\end{figure}

A \textit{Magic Square} is a collection of nine cards, where each card is part of four valid sets~\cite{Falco, Fairbanks} (Fig.~\ref{fig:magicsquare}). Note that for any three non-collinear points (i.e., any three cards that do not form a set), we can construct a magic square, and that each magic square is a \textit{2-flat}, where a $k$-flat is defined by~\cite{Davis} to be a $k$-dimensional affine subspace of a vector space, namely $\F^4$. So a magic square is a \textbf{plane} $p+V$, where $V$ is a two-dimensional subspace and $p$ is a point. A plane has $9$ points, and since two points determine a unique line of $3$ points, it has $\binom{9}{2}/3=12$ lines, in the $4$ directions given by the nonzero vectors of $V$. We use the following standard facts about planes.

\begin{lemma}\label{lem:planes}
Let $P=p+V$ be a plane. (P1)~Through every point of $P$ there is exactly one line of $P$ in each of the four directions of $V$, and these four lines cover all nine points of $P$. (P2)~The line through two points of $P$ lies in $P$. (P3)~Two distinct lines through a common point, or a line and a point not on it, lie in exactly one plane.
\end{lemma}
\begin{proof}
Let $u,v$ be a basis of $V$, so the points of $P$ are $p+au+bv$ with $a,b\in\F$, and $P=q+V$ for every $q\in P$. (P1)~The four directions of $V$ are $u$, $v$, $u+v$, and $u-v$. For each such direction $w$ and each $q\in P$, the line $\{q,q+w,q+2w\}$ lies in $P$, and it is the only line through $q$ with direction $w$. Two of these four lines share only $q$, so together they contain $1+4\cdot 2=9$ points, which is all of $P$. (P2)~If $a=p+x$ and $b=p+y$ with $x,y\in V$, then the third point of their line is $-a-b=-2p-x-y=p-x-y\in P$. (P3)~Every plane through $q$ is $q+W$ for a two-dimensional subspace $W$. Two distinct lines through $q$ have independent directions $u',v'$, and they lie in $q+W$ exactly when $W=\operatorname{span}(u',v')$. A line $\{a,a+w,a+2w\}$ that misses $q$ lies in $q+W$ exactly when $W=\operatorname{span}(a-q,w)$, where $a-q$ and $w$ are independent since $q$ is not on the line. Hence distinct planes through a point $p$ share at most one line through $p$, and a line missing $p$ lies in at most one of them.
\end{proof}

\section{From Four Properties to Three}\label{sec:projection}

In earlier work~\cite{stevens2025maximum}, we proved $M(n,4)=M(n,3)$ for $3\le n\le 12$ by showing a maximum configuration leaves a direction unused. We use the same argument for $13\le n\le 16$.

Consider any configuration $S$ of $n$ cards containing $m=N(S)$ sets. The number of pairs of cards is $\binom{n}{2}$, but each set has three pairs determining the same line. Thus, every set reduces the number $L$ of distinct lines through two points of $S$ by two, giving $L=\binom{n}{2}-2m$. Every direction used by $S$ is the direction of one of these $L$ lines.

\begin{lemma}\label{lem:projection}
If a board $S$ of $n$ cards leaves a direction $v$ unused, then $N(S)\le M(n,3)$.
\end{lemma}
\begin{proof}
Consider the projection $\pi:\F^4\to\F^4/\langle v\rangle\cong\F^3$. Since $v$ is unused, no two selected points differ by a multiple of $v$, and these multiples are exactly the nonzero elements of the kernel of $\pi$. Therefore, the projection preserves all $n$ distinct points. Furthermore, every original set remains a set after projection, since its three images are distinct and $a+b+c=0\implies\pi(a)+\pi(b)+\pi(c)=0$. Since projection preserves distinct points and all original sets, we see that $N(S)\le N(\pi(S))\le M(n,3)$.
\end{proof}

\section{Proofs for $13$, $14$, $15$, and $16$ Cards}\label{sec:main}

We claim that in the cases of $13$, $14$, $15$, and $16$ cards, there will always be an unused direction, apart from boards containing a magic square when $n=15$ or $n=16$; therefore, Lemma~\ref{lem:projection} gives us that $M(n,4)=M(n,3)$ for these values. To find an unused direction, we look for repeated directions: if $k$ of the $L$ lines through two points of $S$ have only $d$ directions, they contain $k-d$ repeated directions, and $S$ uses at most $L-(k-d)$ directions.

\subsection{$n=13$ case}
We know by the computer search of~\cite{stevens2025maximum} that a maximum configuration $S$ contains at least $16$ sets, so $L\le\binom{13}{2}-2\cdot 16$. Since each set contains three cards, there are at least $3\cdot 16=48$ card--set incidences. However, there are only $13$ cards, therefore by the pigeonhole principle, at least one point $p$ must be involved in at least $\lceil \frac{48}{13}\rceil=4$ sets.

If three of these sets lie in a plane that is isomorphic to $\F^2$, then there are $\frac{3^2-1}{2}=4$ distinct directions here, but $7$ points and at most $M(7,4)=5$ sets leading to at least $\binom{7}{2}-2\cdot 5=11$ lines, which lie in the plane by (P2). However, then there are $11-4=7$ duplicate directions. Therefore, our original configuration of 13 cards determines at most $\binom{13}{2}-2\cdot16-7=39<40$ distinct directions.

Alternatively, suppose no three of the four set lines through $p$ lie in a common affine plane. Consider any two of these set lines. They span an affine plane containing five selected points, which determine $\binom{5}{2}-2\cdot 2=6$ distinct lines, since any other set among these five points would contain two points of one set line, and hence equal it. However, this plane has only four possible directions, so there are at least $2$ repeated directions. There are $\binom{4}{2}=6$ pairs of set lines through $p$. Moreover, the four additional lines determined by each pair are distinct from those determined by any other pair. Otherwise, since an additional line avoids $p$ and so lies in only one plane through $p$, three of our set lines would lie in one plane, contrary to our assumption. Thus, these six pairs contribute at least $6\cdot 2=12$ repeated directions, since the four additional lines of each plane add at most two new directions. Therefore, our original configuration of 13 cards determines at most $\binom{13}{2}-2\cdot16-12=34<40$ distinct directions.

In either case, we have shown that fewer than 40 directions are used. Since $\F^4$ has exactly 40 possible directions, there must be at least one unused direction. Due to this, we see that once again $M(13,4)=M(13,3)=16$.

\subsection{$n=14$ case}
Since we see that there are at least $19$ sets for $14$ cards in $\F^4$ by the computer search, a maximum configuration $S$ has at least $19$ sets, and the original counting argument gives us $L\le\binom{14}{2}-2\cdot 19=53$ distinct lines. If we can identify at least $14$ repetitions of directions among them, then this will take us below the $40$ possible directions, and show us once again that at least one direction must be unused.

There are at least $3\cdot 19=57$ card--set incidences. Since there are only $14$ cards, the pigeonhole principle implies that some card $p$ is involved in at least $\lceil \frac{57}{14} \rceil=5$ sets.

Now consider the $\binom{5}{2}=10$ pairs of set lines through $p$. Each pair spans exactly one plane, and as in the $n=13$ case, the repeated directions found in different planes add up: apart from the set lines through $p$, every line we count avoids $p$ and so lies in only one plane through $p$. Thus, if a plane contains $c$ set lines through $p$ and $L_c$ determined lines, its $L_c-c$ remaining lines introduce at most $4-c$ new directions. Hence it contributes at least $L_c-4$ repeated directions.

If each pair spans a distinct plane, then notice that each plane has at most $4$ directions in it, however, each pair contributes $\binom{5}{2}-2\cdot 2=6$ lines, so there are at least $6-4=2$ repeated directions. Since there are $10$ pairs, there are at least $20$ repeated directions in this case. This gives at most $53-20=33$ distinct directions, which is less than $40$, so at least one direction is unused.

If $3$ of the sets lie in the same plane, then there are $7$ points in their union, so as we argued in the $n=13$ case, this leads to at least $11$ lines, but only $4$ possible directions; therefore, this plane contributes $11-4=7$ repeated directions, compared to $\binom{3}{2}\cdot 2=6$ for its three pairs; hence this increases the number of repeated directions. Once again, we have an unused direction.

If $4$ of the sets lie in the same plane (notice that $5$ cannot, since there are only $4$ different directions in the plane), then there are $9$ points. This gives a maximum of $12$ sets, so there are $\binom{9}{2}-2\cdot 12=12$ lines, and thus $12-4=8$ repeated directions. Furthermore, the fifth set and each set on the plane span another plane; hence they contribute another $4\cdot 2=8$ repeated directions, for a total of $16$ repeated directions. This still gives at most $53-16=37<40$ distinct directions.

So in every case there is an unused direction, and we map the projection from $\F^4$ to $\F^3$ and get at most $19$ sets.

\subsection{$n=15$ case}\label{sec:fifteen}

A maximum configuration $S$ has at least $23$ sets, so there are at most $L\le\binom{15}{2}-2\cdot 23=59$ lines, and we must identify at least $20$ repetitions of directions for there to be an unused direction. As in the $n=14$ case, the pigeonhole argument still applies: there are at least $3\cdot 23=69$ card--set incidences, and so at least one card $p$ is involved in at least $\lceil \frac{69}{15} \rceil=5$ sets. 

If no four of the five set lines through $p$ lie in a common plane, then by the $n=14$ count every plane contributes at least $2$ repeated directions per pair, so we can identify at least $20$ repetitions of directions, and the projection reduces this case to the $\F^3$ case again. The same holds when exactly $3$ of the set lines lie in a plane, and we get even more repeated directions as we argued in the $n=14$ case previously.

However, there is one case in which we can only identify $16$, which is when we have a magic square contained within the points: by the $n=14$ count, this is when a plane contains $4$ of the set lines through $p$, and then by (P1) all nine points of this plane are in $S$.

In this case, consider the magic square $G$ plus the $6$ additional points $\{s_1, s_2, \cdots, s_6\}$. The only two new possible sets are $1$ point in the magic square plus $2$ points outside it (we call these \emph{mixed sets}), or entirely within the $6$ new points, since by (P2) a set with two points in $G$ lies in $G$.

Since six cards can form at most $M(6,4)=3$ sets, we can have at most $3$ sets outside of the magic square. Furthermore, translating every card by the same vector maps sets to sets, so we may assume that $G$ is a two-dimensional subspace $V$. Consider the $3^{4-2}=9$ planes parallel to $G$, including $G$. The other $8$ are of the form $G\pm v_i$ for $i\in \{1, 2, 3, 4\}$, which are the cosets. If $s$ and $t$ outside of $G$ form a set with a point $u\in G$, then $s+t=-u\in G$, so $s$ and $t$ belong to opposite cosets, that is, they lie in planes $G+v$ and $G-v$, respectively.

Let $a_i$ and $b_i$ be the number of points in $G+v_i$ and $G-v_i$, then we see that $\sum_{i=1}^{4}(a_i+b_i)=6$ since there are $6$ additional points outside the magic square. However, each such pair lies in at most one set, so the number of mixed sets is at most $\sum_{i=1}^{4}a_ib_i\le \sum_{i=1}^{4}a_i\sum_{i=1}^{4}b_i$, since these are all nonnegative, which is at most $3\cdot 3=9$ by the AM-GM inequality. If there are at most $8$ mixed sets, then we have at most $12+8+3=23$ sets, as desired.

Otherwise, if there are $9$ mixed sets, then both inequalities are equalities, so $\sum_i a_i=\sum_i b_i=3$ and $a_ib_j=0$ for $i\neq j$. Hence $a_i=b_i=3$ for a single $i$, and three of the points are in one coset $G+v$ and the other three are in the other coset $G-v$. Notice that there are at most two sets in this case since all three points must be in the same coset. If two points are in one coset and one is in the other, assume WLOG that two of the points are in $G+v$ and the other is in $G-v$. Then suppose $x,y\in G+v$ and $z\in G-v$. Hence $x=v+a$, $y=v+b$, $z=-v+c$, where $a,b,c$ are all in $V$, and $v\not\in V$. Since we assume they form a set, we get $0=x+y+z=v+a+b+c$. Since $a+b+c\in V$, since $V$ is a vector subspace, we must also have $v=-(a+b+c)\in V$, contradiction. The same argument applies if $x,y\in G-v$ and $z\in G+v$. Hence, there can be at most two additional sets, one in each coset, giving a maximum number of sets again as $12+9+2=23$ sets.

\subsection{$n=16$ case}\label{sec:sixteen}

A maximum configuration $S$ has at least $26$ sets by the computer search. Suppose that $S$ has at least $27$ sets. Then $L\le\binom{16}{2}-2\cdot 27=66$, so we must identify at least $27$ repetitions of directions for there to be an unused direction. There are at least $3\cdot 27=81$ card--set incidences, so some card $p$ is involved in at least $\lceil\frac{81}{16}\rceil=6$ sets. If no four of the six set lines through $p$ lie in a common plane, then by the $n=14$ count the $\binom{6}{2}=15$ pairs give at least $30$ repeated directions, so at most $66-30=36<40$ directions are used, and the projection gives $N(S)\le M(16,3)=26$, a contradiction. Otherwise, as in the $n=15$ case, four of these set lines lie in a plane, and by (P1) $S$ contains this magic square $G$.

So $S$ consists of a magic square $G$ and $7$ further cards. As in the $n=15$ case, we may assume that $G=V$, and every set of $S$ is one of the $12$ sets of $G$, a mixed set with one card in $G$ and its other two cards in opposite cosets, or a set among the $7$ cards outside $G$. Let $a_i$ and $b_i$ be the numbers of these cards in $G+v_i$ and $G-v_i$, so $\sum_{i=1}^{4}(a_i+b_i)=7$. As in the $n=15$ case, there are at most $\sum_{i=1}^4 a_ib_i$ mixed sets.

If the $7$ cards meet at least two of the four pairs of cosets, let $s_i=a_i+b_i$, and let $s$ be the largest $s_i$, so $1\le s\le 6$. Since $a_ib_i\le s_i^2/4$ by AM-GM, and the other $s_i$ sum to $7-s$, there are at most $\sum_i a_ib_i\le\frac14\bigl(s^2+(7-s)^2\bigr)\le\frac14(6^2+1^2)<10$ mixed sets. Since $7$ cards form at most $M(7,4)=5$ sets, $N(S)\le 12+9+5=26$, a contradiction.

Otherwise, all $7$ cards lie in one pair of cosets $G+v$ and $G-v$, with $a$ and $b=7-a$ cards. There are at most $ab$ mixed sets. By the argument of the $n=15$ case, no set has two cards in one of these cosets and its third card in the other, so every other set lies in $G+v$ or in $G-v$, and there are at most $M(a,4)+M(b,4)$ of them. Hence $N(S)\le 12+ab+M(a,4)+M(b,4)$, and up to swapping $a$ and $b$, the splits $7+0$, $6+1$, $5+2$, and $4+3$ give at most $12+0+5=17$, $12+6+3=21$, $12+10+2=24$, and $12+12+2=26$ sets. In every case $N(S)\le 26$, which contradicts $N(S)\ge 27$. Hence $M(16,4)=M(16,3)=26$, and together with~\cite{stevens2025maximum}, this proves the Main Theorem.

\subsection{Constructions}

\begin{figure}[t]
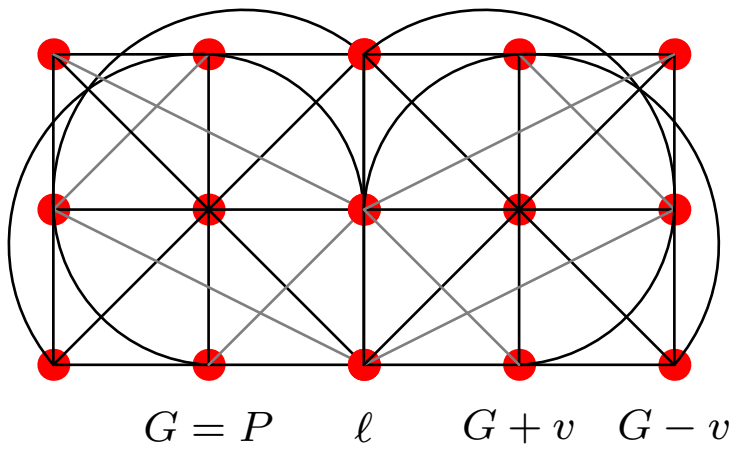

\centering
  \begin{asy}
    import graph;
    size(1.9inch);
    defaultpen(fontsize(8pt));
    picture sq;
    dot(sq, (0,0), 6+red);
    dot(sq, (1,0), 6+red);
    dot(sq, (0, 1), 6+red);
    dot(sq, (0,2), 6+red);
    dot(sq, (1, 1), 6+red);
    dot(sq, (1,2), 6+red);
    dot(sq, (2,1), 6+red);
    dot(sq, (2,0), 6+red);
    dot(sq, (2,2), 6+red);
    draw(sq, (0,0)--(0,1)--(0,2), black);
    draw(sq, (1,0)--(1,1)--(1,2), black);
    draw(sq, (0,1)--(1,1)--(2,1), black);
    draw(sq, (0,0)--(1,0)--(2,0), black);
    draw(sq, (0,2)--(1,1)--(2,0), black);
    draw(sq, (2,0)..(1,2)..(0,1), black);
    draw(sq, (0,0)--(2,1)--(1,2), grey);
    draw(sq, (1,0)..(2,1)..(0,2), black);
    draw(sq, (2,0)--(2,1)--(2,2), black);
    draw(sq, (0,2)--(1,2)--(2,2), black);
    draw(sq, (0,0)--(1,1)--(2,2), black);
    draw(sq, (1,0)--(0,1)--(2,2), grey);
    add(reflect((1,0), (1,1))*sq);
    add(shift(2,0)*sq);
    label("$G=P$", (1,-0.4));
    label("$\ell$", (2,-0.4));
    label("$G+v$", (3,-0.4));
    label("$G-v$", (4,-0.4));
  \end{asy}
  \caption{The board $P\cup Q$: the magic squares $P$ (left) and $Q$ (right) share the middle column $\ell$, so there are $12+12-1=23$ sets. The two right columns are the sets within $G\pm v$, and the $9$ non-vertical sets of $Q$ are the mixed sets.}
  \label{fig:pq}
\end{figure}

For $x,y,z\in\F$, let $P=\{(x,y,0,0)\}$ and $Q=\{(x,0,z,0)\}$ (Fig.~\ref{fig:pq}). These planes share the line $\ell=\{(x,0,0,0)\}$, so $P\cup Q$ has $15$ cards. A line contained in $P\cup Q$ has two points in $P$ or two in $Q$, so by (P2) it lies in $P$ or in $Q$. Hence the sets of a board inside $P\cup Q$ are the lines of $P$ and of $Q$ that lie in the board, and $P\cup Q$ has $12+12-1=23$ sets.

For $15$ cards, $P\cup Q$ is the case of $9$ mixed sets in Section~\ref{sec:fifteen}: with $G=P$ and $v=(0,0,1,0)$, we have $G+v=\{(x,y,1,0)\}$ and $G-v=\{(x,y,2,0)\}$, so the three cards $(x,0,1,0)$ lie in $G+v$ and the three cards $(x,0,2,0)$ lie in $G-v$. The $9$ pairs $(x,0,1,0)$, $(x',0,2,0)$ form sets with $(-x-x',0,0,0)\in\ell$, and each of the two triples forms a set, so $12+9+2=23$.

For $14$ and $13$ cards, we remove cards of $P\setminus\ell$. By (P1), $(0,1,0,0)$ lies on exactly $4$ sets, all in $P$, so removing it leaves $14$ cards with $23-4=19$ sets. Also removing $(0,2,0,0)$ removes $4$ more sets, one of which, $\{(0,y,0,0)\}$, is already gone, so $13$ cards remain with $19-3=16$ sets.

For $16$ cards, we add $w=(0,1,1,0)\in G+v$, the $4+3$ case of Section~\ref{sec:sixteen}. The card $w$ forms three mixed sets $\{w,(x,0,2,0),(-x,2,0,0)\}$, one with each card of $G-v$. No additional set lies entirely within $G+v$, since the cards $(x,0,1,0)$ have second coordinate $0$, whereas $w$ has second coordinate $1$. Thus there are $23+3=26$ sets.

\section{Conclusion and Open Questions}

It remains an open question to find, with proof, the maximum number of sets with $4$ properties for $n>16$. We conjecture that the values of $M(n,3)$ in OEIS A380167~\cite{oeis_A380167} match the maximum for $4$ properties. For $n=17$, where $M(17,3)=30$, our argument falls short: a board with $31$ sets determines at most $\binom{17}{2}-2\cdot 31=74$ lines, so an unused direction needs $35$ repeated directions, but a card on $\lceil 93/17\rceil=6$ sets gives only $30$. Finding the maximum number of sets with more than $4$ properties also remains open.

\section*{Acknowledgements}
We thank Duncan Wilson and Richard Kelley, coauthors of the conference paper on which this work builds; Haydn Bradstreet and the anonymous reviewers of the conference paper for their comments; and the OEIS team for discussing the sequence. The AI assistants ChatGPT (OpenAI) and Claude (Anthropic) helped develop the proofs, including the argument for $16$ cards, write and edit the text, and draw Fig.~\ref{fig:pq}; the author wrote the original draft and verified all the content.

\bibliographystyle{IEEEtran}
\bibliography{sample}

\end{document}